\documentclass[preprint,12pt,authoryear]{elsarticle}
\usepackage{amssymb}
\RequirePackage{fix-cm}
\usepackage{graphicx}
\usepackage{color}
\usepackage{caption}
\usepackage{amsmath}
\usepackage{breakcites}
\usepackage{multicol} 
\usepackage{url}
\usepackage{amsthm}
\usepackage{comment}
\usepackage{natbib}
\usepackage{mathrsfs}

\numberwithin{equation}{section}
\theoremstyle{plain}

\newtheorem*{prop47}{Theorem 1}
\newtheorem*{prop48}{Theorem 2}
\newtheorem*{prop46}{Lemma 2}
\newtheorem*{prop45}{Lemma 1}
\newcommand*\cpp{C\kern-0.2ex\raisebox{0.4ex}{\scalebox{0.8}{+\kern-0.4ex+}}}
\newcommand{\slfrac}[2]{\left.#1\middle/#2\right.}

\DeclareMathOperator*{\cov}{cov}
\DeclareMathOperator*{\ave}{ave}

\DeclareMathOperator*{\argmin}{\arg\!\min}
\DeclareMathOperator*{\argmax}{\arg\!\max}
\journal{Statistics \& Probability Letters}
\begin{document}

\begin{frontmatter}
\title{The Finite Sample Breakdown Point of PCS}
\author[eric]{Eric Schmitt}
\ead{eric.schmitt@wis.kuleuven.be}
\author[vik]{Viktoria \"Ollerer}
\ead{viktoria.oellerer@kuleuven.be}
\author[kav]{Kaveh Vakili}
\ead{kaveh.vakili@wis.kuleuven.be}

\address[eric]{\emph{Corresponding author}. Afdeling Statistiek, Celestijnenlaan 200b - bus 2400, 3001 Leuven. Tel +32 16 37 23 40.}
\address[vik]{Faculty of Business and Economics, ORSTAT, KU Leuven, Belgium.}
\address[kav]{Afdeling Statistiek, Celestijnenlaan 200b - bus 2400, 3001 Leuven.}

\begin{abstract}
The Projection Congruent Subset (PCS) is new method for finding multivariate outliers. PCS 
 returns an outlyingness index which can be used to construct affine equivariant estimates of 
 multivariate location and scatter. In this note, we derive the finite sample breakdown point 
 of these estimators.

\end{abstract}
\begin{keyword}
breakdown point \sep robust estimation \sep multivariate statistics.
\end{keyword}
\end{frontmatter}

\section{Introduction}
Outliers are observations that depart from
the pattern of the majority of 
the data. 
Identifying outliers is a major concern 
in data analysis because a few outliers, 
if left unchecked, can exert a 
 disproportionate pull on the fitted 
parameters of any statistical model, 
preventing the analyst from uncovering
the main structure in the data. 

To measure the robustness of an estimator to
the presence of outliers in the data,~\cite{ppcs:D82} 
introduced the notion of finite sample breakdown point. 
Given a sample and an estimator, this
is the smallest number of observations 
 that need to be replaced by outliers to cause the 
 fit to be arbitrarily far from the values it would
 have had on the original sample.
Remarkably, the finite sample breakdown point
of an estimator can be derived 
without recourse to concepts of chance or randomness
using geometrical features of a sample alone~\citep{ppcs:D82}.
Recently,~\cite{hcs:VS13} introduced 
the Projection Congruent Subset (PCS) method. 
PCS computes an outlyingness index, as well as 
estimates of location and scatter derived from it. 
The objective of this paper is to establish 
the finite sample breakdown of these estimators
and show that they are maximal.

Formally, we begin from
the situation whereby the data matrix $\pmb X$, 
is a collection of  
$n$ so called $genuine$ observations 
drawn from a $p$-variate model $F$ with $p>1$. 
However, we do not observe $\pmb X$ 
but an $n \times p$ (potentially) corrupted data set 
$\pmb X^{\varepsilon}$ that consists of
$g<n$ observations from $\pmb X$ and 
$c=n-g$  arbitrary values,
with $\varepsilon=\slfrac{c}{n}$, 
denoting the (unknown) rate of 
contamination. 

Historically, the goal of many robust estimators has been to achieve high breakdown while obtaining reasonable efficiency.
PCS belongs to a small group of robust estimators that 
have been designed to also have low bias (see~\cite{ppcs:MA92},~\cite{ppcs:AY02} and~\cite{ppcs:AY10}). 
In the context of robust estimation, a 
low bias estimator reliably finds a fit 
close to the one it would 
have found without the outliers, when $c\leqslant n-h$
with $h=\lceil(n+p+1)/2\rceil$. 
To the best of our knowledge, PCS is  
the first member of this group of estimators
to be supported by a fast and affine equivariant 
algorithm (FastPCS) enabling its use by practitioners.

The rest of this paper unfolds as follows. 
In Section~\ref{s2}, we detail the PCS 
estimator. In Section~\ref{s3}, we formally detail the concept 
of finite sample breakdown point of an estimator 
and establish the notational conventions we will 
use throughout. Finally, in Section~\ref{s4}, we prove  
the finite sample breakdown point of PCS.

\section{The PCS criterion}\label{s2}

Consider a potentially contaminated data set $\pmb X$ 
 of $n$ vectors 
$\pmb x_i\in\mathbb{R}^p$,  with $n>p+1>2$.
Given all $M=\binom{n}{h}$ possible 
$h$-subsets $\{H^m\}_{m=1}^M$,
PCS looks for the one   
 that is most {\em congruent} 
along many univariate projections. 
Formally, given an $h$-subset $H^m$, we 
denote $B(H^m)$ the set of all vectors normal
to hyperplanes spanning a $p$-subset of $H^m$.
More precisely, all directions $\pmb a \in B(H^m)$ define hyperplanes $\{\pmb x\in\mathbb{R}^p: \pmb x'\pmb a=1\}$ that contain $p$ observations of $H^m$. For $\pmb a \in B(H^m)$ and $\pmb x_i\in\pmb X$, we can compute the squared orthogonal distance, $d_i^2$, of $\pmb x_i$ to the hyperplane defined by $\pmb a$ as
\begin{equation}\label{pcs:crit0}
    d_i^2(\pmb a) =\frac{(\pmb a'\pmb x_i-1)^2}{||\pmb a||^2}\;.
\end{equation}
The set of the $h$ observations with smallest $d_{i}^2(\pmb a)$ is then defined as
\begin{equation}
H^{\pmb a}=\{i:d_{i}^2(\pmb a)\leqslant d_{(h)}^2(\pmb a)\},
\end{equation}
where $d_{(h)}$ denotes the $h$th-order statistic of a vector $\pmb d$.

We begin by considering the case in which 
 $d_{(h)}^2(\pmb a)>0$. For a given subset $H^m$ and direction $\pmb a$ we define the {\em incongruence index} 
of $H^m$ along $\pmb a$ as
\begin{equation}\label{pcs:crit1}
   I(H^m,\pmb a):=
    \log\left(\frac{\displaystyle\ave_{i\in H^m}d^2_i(\pmb a)}{\displaystyle\ave_{i\in H^{\pmb a}}d_i^2(\pmb a)}\right)\;
\end{equation}
with the conventions that $\log(0/0):=0$. This index is always positive and will be smaller the more members of $H^m$ correspond with, or are similar to, the members of $H^{\pmb a}$.
To remove the dependency of Equation \eqref{pcs:crit1} on $\pmb a$, we measure the incongruence of $H^m$ by considering the average over many directions $\pmb a\in B(H^m)$ as
\begin{equation}\label{mpcs}
      I(H^m):=\ave_{\pmb a\in B(H^m)} I(H^m,\pmb a).\;
\end{equation}
The optimal $h$-subset, $H^*$, is the one satisfying the PCS criterion:
\begin{equation}
H^*=\underset{\{H^m\}_{m=1}^M}{\argmin}\;I(H^m).\nonumber
\end{equation}
Then, the {\em PCS estimators of location and scatter} are the sample mean and covariance of the observations with indexes in $H^*$:
 \begin{equation}
\left(\pmb t^*(\pmb X),\pmb S^*(\pmb X)\right)=\left(\ave_{i\in H^*}\pmb x_i,\cov_{i\in H^*}\pmb x_i\right). \nonumber
\end{equation}

Finally, we have to account for the special case where $d_{(h)}^2(\pmb a)=0$. In this case, we enlarge $H^*$ to be the subset of all observations lying on $\pmb a$. More precisely, if $\exists\;\pmb a^*\in B(H^*):|\{i:d_i^2(\pmb a^*)=0\}|\geqslant h$, then $H^*=\{i:d_i^2(\pmb a^*)=0\}$.\\

\subsection{Illustrative Example} 

\begin{figure}[!ht]
   \centering
    \includegraphics[width=0.9\textwidth]{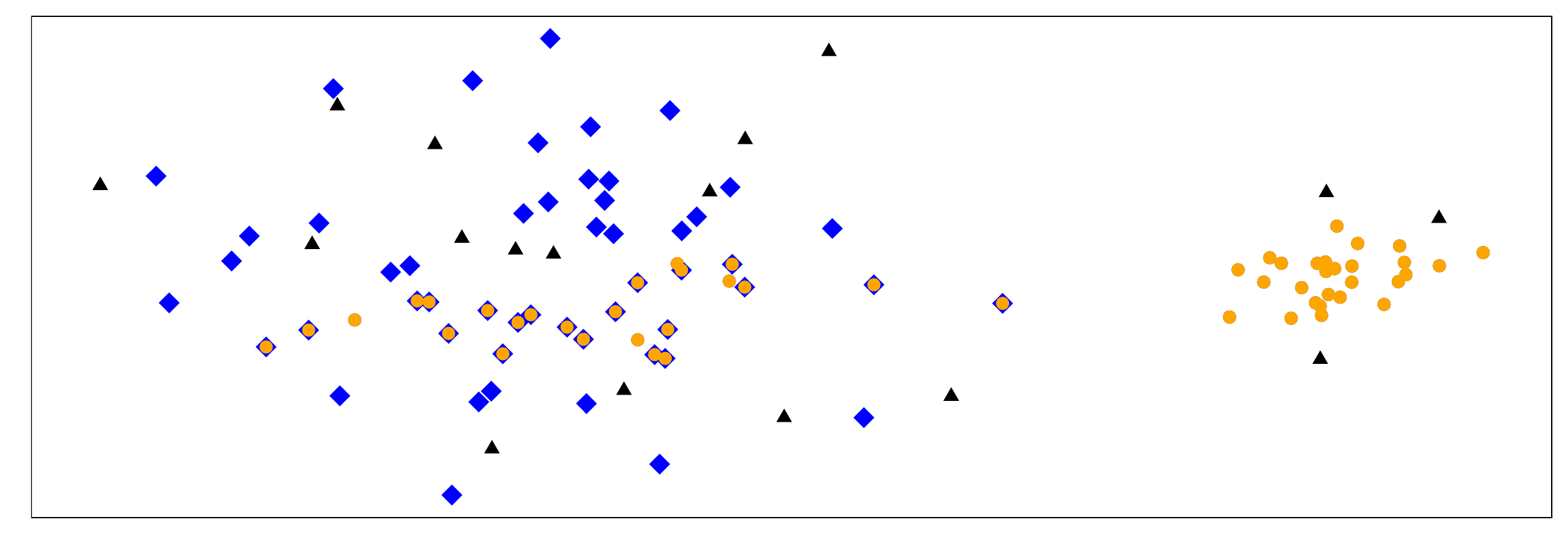}
  \caption{Bivariate data example. The members of $H^1$ ($H^2$) are depicted as dark blue diamonds (light orange circles).}
    \label{fig:dataPlot}
\end{figure}

To give additional insight into PCS and the characterization of a cloud of point in terms of congruence, we provide the following example. Figure~\ref{fig:dataPlot} depicts a data set $\pmb X^\varepsilon$ of 100 observations, 30 of which come from a cluster of outliers on the right. 
For this data set, we draw two $h$-subsets of 52 observations each. 

Subset $H^1$ (dark blue diamonds) contains only genuine observations, while subset $H^2$ (light orange circles) contains 27 outliers and 25 genuine observations. Finally, 
the 17 observations belonging to neither $h$-subset are depicted as black triangles.
For illustration's sake, we selected the members of $H^2$ 
 so that their covariance has smaller determinant than $any$ $h$-subsets formed of genuine 
 observations. 
Consequently, robust methods based on a characterization of $h$-subsets in terms of density
alone will always prefer the contaminated subset $H^2$ over any uncontaminated $h$-subset 
(and in particular $H^1$).

The outlyigness index computed by PCS differs from that of other robust estimators in two important ways. 
First, in PCS, the data is projected onto directions given by $p$ points
 drawn from the members of a given subset, $H^m$, rather than indiscriminately from the entire data set. This choice is motivated by the fact that when $\varepsilon$ and/or $p$ are high, the vast majority of random $p$-subsets of $\{1,\dots,n\}$ will be contaminated. If the outliers are concentrated, this yields directions almost parallel to each other. In contrast, for an uncontaminated $H^m$, our sampling strategy always ensures a wider spread of directions and this yields better results.

\begin{figure}[!ht]
   \centering
    \includegraphics[width=0.9\textwidth]{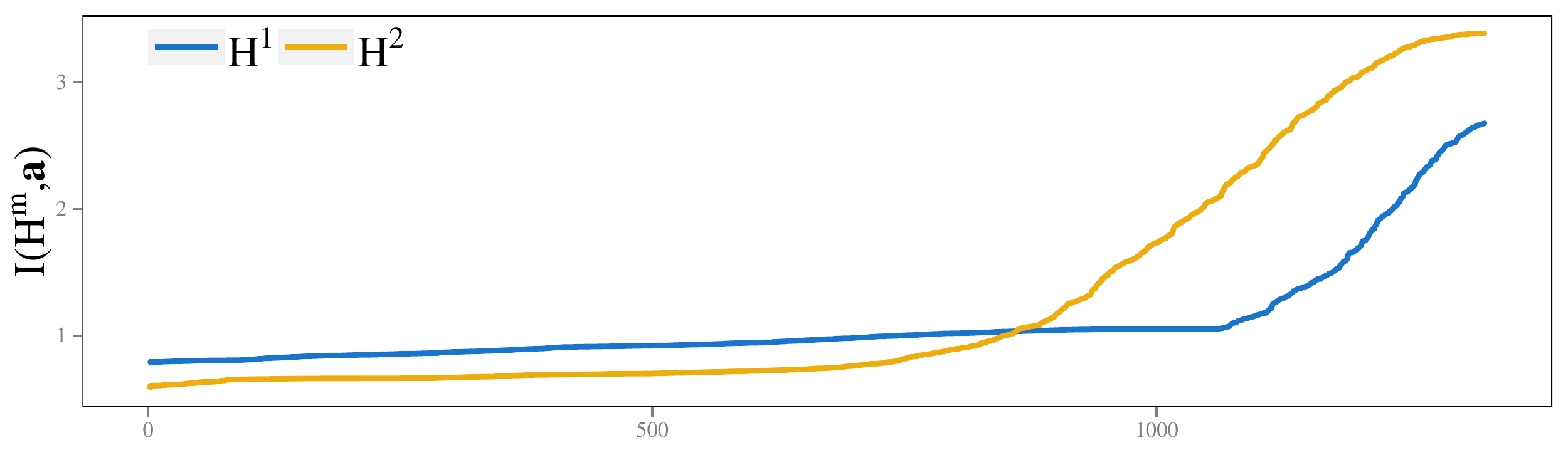}
  \caption{Sorted values of $I(H^m,\pmb a)$ for $H^1$ and $H^2$, shown as dark blue (light orange) lines.}
    \label{fig:Icompare}
\end{figure}

The second feature of PCS is that the congruence index 
used to characterize an $h$-subset depends on all the data points in the sample. 
We will illustrate this by considering all $\binom{52}{2}=1326$ members of $B(H^1)$. For 
 each, we compute the corresponding value of $I(H^1, \pmb a)$. Then, we sort these and plot them in Figure~\ref{fig:Icompare}. We do the same for $H^2$. We note in passing that $I(H^1)<I(H^2)$.

Consider now in particular the values of the $I$-index corresponding to $H^1$ and starting at around 1050 on the horizontal 
axis of Figure~\ref{fig:Icompare}. 
These higher values of $I(H^1,\pmb a)$ correspond to members of $B(H^1)$ that are aligned with the vertical axis (i.e. they correspond to horizontal hyperplanes), and are much larger than the remaining values of $I(H^1,\pmb a)$.  
This is because, for the data configuration shown in Figure~\ref{fig:dataPlot}, the outliers 
do not stand out from the good data in terms of their orthogonal distances to hyperplanes defined by the vertical directions. 
  As a result, this causes many outliers to enter the sets $H^{\pmb a}$ and this deflates the values of the 
   $\displaystyle\ave_{i\in H^{\pmb a}}d_i^2(\pmb a)$ corresponding to these directions. 
   Since the set $H^1$ is fixed, there is no corresponding effect on $\displaystyle\ave_{i\in H^1}d^2_i(\pmb a)$ so that the outliers will  influence the values of $I(H^1,\pmb a)$ for some directions $\pmb a$,
    even though $H^1$ itself is uncontaminated. 
This apparent weakness is an inherent feature of PCS. 
In the remainder of this note, we prove the following 
 counter-intuitive fact: outliers influence the value of $I(H^m)$, even when $H^m$ is free of outliers,
 yet, so long as there are fewer than $n-h$ of them, their influence on the PCS fit will always remain bounded. In other words, breakdown only occurs if $c>n-h$ (see Section~\ref{s4}).

\section{Finite sample breakdown point}\label{s3}
To lighten notation and without loss of generality, we arrange the observed data matrix $\pmb X^{\varepsilon}=((\pmb X^g)', (\pmb X^{c})')'$ 
with rows $\{\pmb x_i^{\varepsilon}\}_{i=1}^n$ so that the $\varepsilon\%$ of contaminated observations $\pmb X^c$ are in the last $c$ rows and the uncontaminated observations $\pmb X^g$ in the first $g$ rows. Then,  $\mathcal{X}^{\varepsilon}$ will refer to the set of all corrupted data sets $\pmb X^{\varepsilon}$ 
and $\mathcal{H}$ is the set of all $h$-subsets of $\{1,\ldots,n\}$, $\mathcal{H}^c=\{H\in\mathcal{H}:H\cap \{g+1,\ldots,n\}\neq\varnothing\}$  the set of all $h$-subsets of $\{1,\ldots,n\}$ with at least one contaminated observation,  and $\mathcal{H}^g=\{H\in\mathcal{H}:H\cap \{g+1,\ldots,n\}=\varnothing\}$ the set of all uncontaminated $h$-subsets of $\{1,\ldots,n\}$.

The following assumptions (as per, for example~\cite{ppcs:T94}) all pertain to the original, uncontaminated,  
data set $\pmb X$. 
In the first part of this note, we will consider the case whereby the point cloud formed by $\pmb X^g=\{\pmb x_i^g\}_{i=1}^g$ lies in \textit{general position} in 
$\mathbb{R}^p$. The following definition of \textit{general position} is adapted from 
\cite{mcs:RL87}:

 \bigskip

\textsc{Definition} 1: \textit{General position in $\mathbb{R}^p$}. 
$\pmb X$ is in general position in $\mathbb{R}^p$ if 
no more than $p$-points of $\pmb X$
 lie in any $(p-1)$-dimensional 
 affine subspace. 
For $p$-dimensional data, this
means that there are no more than $p$ points
 of $\pmb X$ on any hyperplane, so that any $p+1$ points of
$\pmb X$ always determine a $p$-simplex with non-zero determinant.

\bigskip
\noindent Throughout, we will also assume that
\begin{eqnarray*}
\underset{\pmb a\in\mathbb{R}^p_{\neq 0}}{\sup}\max\left(\frac{\pmb X'\pmb a}{||\pmb a||}\right)^2\text{ is bounded}
\end{eqnarray*}
and that the genuine observations contain no duplicates:
\begin{eqnarray*}
||\pmb x_i-\pmb x_j||>0\forall\;1\leqslant i<j\leqslant n. 
\end{eqnarray*}
For any $h$-subset $H^m\in\mathcal{H}$ and 
$\pmb X^{\varepsilon}$, we will denote the sample mean and covariance of the observations with indexes in $H^m$ as
 \begin{equation}
\left(\pmb t^m(\pmb X^{\varepsilon}),\pmb S^m(\pmb X^{\varepsilon})\right)=\left(\ave_{i\in H^m}\pmb x_i^{\varepsilon},\cov_{i\in H^m}\pmb x_i^{\varepsilon}\right). \nonumber
\end{equation} 
Then, given $\pmb X^{\varepsilon}\in\mathcal{X}^{\varepsilon}$, and an affine equivariant estimator of location $\pmb t$, we define the bias of $\pmb t$ at $\pmb X$ as
\begin{eqnarray}
\text{bias}(\pmb t,\pmb X,\varepsilon)=\underset{\pmb X^\varepsilon\in\mathcal{X}^{\varepsilon}}{\sup}\;||\pmb t(\pmb X^{\varepsilon}) - \pmb t(\pmb X)||. \nonumber
\end{eqnarray}
Furthermore, given $\pmb X^{\varepsilon}\in\mathcal{X}^{\varepsilon}$, genuine data $\pmb X$ and an affine equivariant estimator of scatter $\pmb S$ with $\pmb S(\pmb X)$ positive definite (denoted from now on by $\pmb S(\pmb X)\succ 0$), we define the bias of $\pmb S$ at $\pmb X$ as
\begin{eqnarray}
\text{bias}(\pmb S,\pmb X,\varepsilon)=\underset{\pmb X^\varepsilon\in\mathcal{X}^{\varepsilon}}{\sup}\;
\frac{\lambda_1(\pmb Q^{\varepsilon})}{\lambda_p(\pmb Q^{\varepsilon})},\nonumber
\end{eqnarray}
where $\pmb Q^{\varepsilon}=(\pmb S(\pmb X)^{-1/2}\pmb S(\pmb X^{\varepsilon})\pmb S(\pmb X)^{-1/2})\succ0$ and 
$\lambda_1(\pmb Q^{\varepsilon})$ ($\lambda_p(\pmb Q^{\varepsilon})$) denotes the largest (smallest) 
eigenvalue of a matrix $\pmb Q^{\varepsilon}$. 
Since PCS is affine equivariant (see Appendix 1), w.l.o.g., we can set $\pmb t(\pmb X) = \pmb 0$ so that the expression of bias reduces to
\begin{eqnarray}
\text{bias}(\pmb t,\pmb X,\varepsilon)=\underset{\pmb X^\varepsilon\in\mathcal{X}^{\varepsilon}}{\sup}\;||\pmb t(\pmb X^{\varepsilon})||. \nonumber
\end{eqnarray}
Furthermore, if that data is in general position and 
$\pmb S(\pmb X)$ is affine equivariant then we can  
w.l.o.g. set $\pmb S(\pmb X)=\pmb I_p$ 
($\pmb I_p$ is the rank $p$ identity matrix) so that the 
expression of the bias reduces to
\begin{eqnarray}
\text{bias}(\pmb S,\pmb X,\varepsilon)= \underset{\pmb X^\varepsilon\in\mathcal{X}^{\varepsilon}}{\sup}\;\frac{\lambda_1(\pmb S(\pmb X^{\varepsilon}))}{\lambda_p(\pmb S(\pmb X^{\varepsilon}))}. \nonumber
\end{eqnarray}
The finite sample breakdown points $\varepsilon^*_n$~\citep{pcs:MMY06} 
 of $\pmb t$ and $\pmb S$ are then defined as
\begin{eqnarray*}
\varepsilon^*_n(\pmb t,\pmb X)&=&\underset{1\leqslant c\leqslant n}{\min}\left\{\varepsilon=\frac{c}{n}:\text{bias}(\pmb t,\pmb X,\varepsilon)=\infty\right\} \\
\varepsilon^*_n(\pmb S,\pmb X)&=&\underset{1\leqslant c\leqslant n}{\min}\left\{\varepsilon=\frac{c}{n}: \text{bias}(\pmb S,\pmb X,\varepsilon)=\infty\right\}.
\label{eq:breakdownexp}
\end{eqnarray*}
Finally, for point clouds $\pmb X$ lying in general position in $\mathbb{R}^p$,~\cite{ppcs:D87} gives a strict upper bound for the finite sample breakdown point for any affine equivariant location and scatter statistics $(\pmb t,\pmb S)$, namely:
 \begin{eqnarray}
 \label{eq:maxBP}
 \varepsilon^*_n(\pmb t,\pmb X)&\leqslant&(n-h+1)/n\notag\\
 \varepsilon^*_n(\pmb S,\pmb X)&\leqslant&(n-h+1)/n
 \end{eqnarray}

\section{Finite sample breakdown point of PCS}\label{s4}

To establish the breakdown point of $\pmb S^*$, we first introduce two lemmas describing properties of the $I$-index.
Both deal with the case where $\pmb X$ lies in general position 
 in $\mathbb{R}^p$. 
Then, we discuss the case where $\pmb X$ does not lie in 
general position. 

In the first lemma, we show that the incongruence index of a clean $h$-subset is bounded. 

\begin{prop45}
Let $c\leqslant n-h$ and $\pmb X$ lies in general position in $\mathbb{R}^p$. Then
\begin{eqnarray}\label{pcs:l1_state}
\underset{\pmb X^{\varepsilon}\in\mathcal{X}^{\varepsilon}}{\sup}\;\underset{H^g\in\mathcal{H}^g}{\max}\;\underset{\pmb a\in B(H^g)}{\max}I(H^g,\pmb a)\leqslant k(\pmb X)
\end{eqnarray}
for any fixed, positive scalar $k(\pmb{X})$ not 
depending on the outliers. 
\end{prop45} 

\begin{proof}
Consider first the numerator of $I(H^g,\pmb a)$. For a fixed $H^g\in\mathcal{H}^g$, we can find for each $\pmb a\in B(H^g)$, the $p$ observations of $\pmb X^g$ that lie furthest away from the hyperplane defined by $\pmb a$. The average of their distances (as given by Equation \eqref{pcs:crit0}) to the hyperplane $\pmb a$ is finite and constitutes an upper bound on the average distance of any $p$ observations of $\pmb X^g$ to the hyperplane $\pmb a$. As we have at most $\binom{h}{p}$ different directions $\pmb a\in B(H^g)$ and only $\binom{n-c}{h}$ uncontaminated subsets $H^g\in\mathcal{H}^g$, the upper bound of the average distances stays finite
\begin{eqnarray}
\underset{H^g\in\mathcal{H}^g}{\max}\;\underset{\pmb a\in B(H^g)}{\max}\ave_{i\in H^g}\;d_i^2(\pmb a)\leqslant U(\pmb{X})\nonumber
\end{eqnarray}
for any positive, fixed, finite scalar $U(\pmb{X})$ not depending on the outliers. Since the contaminated observations have no influence on the distance $d_i^2(\pmb a)$ with $i \in H^g(\pmb a)$ for $\pmb a\in B(H^g)$, we can say that
\begin{eqnarray}\label{pcs:l1_aux1}
\underset{\pmb X^{\varepsilon}\in\mathcal{X}^{\varepsilon}}{\sup}\;\underset{H^g\in\mathcal{H}^g}{\max}\;\underset{\pmb a\in B(H^g)}{\max}\ave_{i\in H^g}\;d_i^2(\pmb a)\leqslant U(\pmb{X}).
\end{eqnarray}

Consider now the denominator of $I(H^g,\pmb a)$. For any $H^g\in\mathcal{H}^g$ and $\pmb a\in B(H^g)$, let $H^{\pmb a}$ denote the subset that consists of the indexes of the $h$ observations of the observed data matrix $\pmb X^\varepsilon$ that lie closest to the hyperplane spanned by $\pmb a$. As $c\leqslant n-h$ and $h=\lceil(n+p+1)/2\rceil$, $H^{\pmb a}$ contains at least $p+1$ uncontaminated observations. In total, when $H^g\in\mathcal{H}^g$ is not fixed, there are at most $\binom{n-c}{p}$ different directions $\pmb a$ defined by a $H^g\in\mathcal{H}^g$. For any $\pmb a$, the smallest value of $\ave_{i\in H^{\pmb a}} d_i^2(\pmb a)$ is attained if the contaminated observations of $H^{\pmb a}$ achieve $d_i^2(\pmb a)=0$. As the uncontaminated observations lie in general position, we know that the $p+1$ uncontaminated observations in $H^{\pmb a}$ cannot lie within the same $p$-dimensional subspace, i.e.
\begin{eqnarray}
\exists i\in H^{\pmb a}: d_i^2(\pmb a)>0. \nonumber
\end{eqnarray}
As the number of uncontaminated observations is fixed, we have that 
\begin{eqnarray}\label{pcs:l1_aux2}
\underset{H^g\in\mathcal{H}^g}{\min}\;\underset{\pmb a\in B(H^g)}{\min}\ave_{i\in H^{\pmb a}}\;d_i^2(\pmb a)\geqslant l(\pmb{X})>0
\end{eqnarray}
for any fixed positive scalar $l(\pmb X)$ not depending on the outliers. This inequality holds even if the outliers have the smallest average distance that is possible (i.e. when $\pmb a:d_i^2(\pmb a)=0$ for the contaminated observations). Thus, Inequality (\ref{pcs:l1_aux2}) holds for any $\varepsilon$-contaminated data set 
$\pmb X^{\varepsilon}$ yielding
\begin{eqnarray}\label{pcs:l1_aux3}\small
\underset{\pmb X^{\varepsilon}\in\mathcal{X}^{\varepsilon}}{\inf}\;\underset{H^g\in\mathcal{H}^g}{\min}\;\underset{\pmb a\in B(H^g)}{\min}\ave_{i\in H^{\pmb a}}\;d_i^2(\pmb a)\geqslant l(\pmb{X})>0.
\end{eqnarray}
Using Equation~\eqref{pcs:crit1} and the Inequalities~\eqref{pcs:l1_aux1} and~\eqref{pcs:l1_aux3}, we get~\eqref{pcs:l1_state}. 
\end{proof}

The second lemma shows the unboundedness of the incongruence index of contaminated subsets.

\begin{prop46}
Let $c\leqslant n-h$ and assume that $\pmb X$ lies in general position in $\mathbb{R}^p$. 
Take a fixed h-subset $H^c\in \mathcal{H}^c$. Then 
\begin{eqnarray}
\forall U_1 >0: \; \exists \pmb X^\varepsilon \in \mathcal{X}^\varepsilon: \;I(H^c,\pmb a)>U_1
\label{eq:l2_aux1}
\end{eqnarray}
for at least one $\pmb a\in B(H^c)$. In other words, for a given set of indexes $H^c$, there exists a data set $\pmb X^\varepsilon$ with contaminated observations with indexes in $H^c$ such that $I(H^c,\pmb a)$ is unbounded.
\end{prop46} 

\begin{proof}
Consider first the numerator of $I(H^c,\pmb a)$. For a fixed $H^c\in\mathcal{H}^c$, denote $G^+=\{G\cap H^c\}$. Since $c\leqslant n-h$, as already mentioned in Lemma 1 above, any $h$-subset contains at least $p+1$ uncontaminated observations, i.e. $|G^+|\geqslant p+1$. Let $B^+(H^c)\subseteq B(H^c)$ be the set of all directions defining a hyperplane spanned 
by a $p$-subset of $G^+$. $|G^+|\geqslant p+1$ yields $|B^+(H^c)|\geqslant p+1$. As the uncontaminated observations $G\supseteq G^+$ lie in general position, the members of 
$B^+(H^c)$ are, by definition, linearly independent. As a result, the outliers can belong to (at most) the subspace spanned by $p$ uncontaminated observations. 
Hence, for every $U_2>0$, there exists at least one member $\pmb a^c_+$ of $B^+(H^c)$, at least one $i\in H^c$ and at least one $\pmb X^{\varepsilon}\in\mathcal{X}^{\varepsilon}$ such that
\begin{eqnarray}
\;d_i^2(\pmb a^c_+)>U_2.
\label{eq:l2_aux2}
\end{eqnarray}

Consider now the denominator of $I(H^c,\pmb a)$: Since the members of $B^+(H^c)$ all pass through  members of $\pmb X^g$ only, we have that
\begin{eqnarray}
d_{(h)}^2(\pmb a^c_+)\leqslant U_3\leqslant\underset{i\leqslant h}{\max}\;d_i^2(\pmb a^c_+).
\label{eq:l2_aux3}
\end{eqnarray}
Using Equation \eqref{pcs:crit1}, and Inequalities \eqref{eq:l2_aux2} and \eqref{eq:l2_aux3}, we get \eqref{eq:l2_aux1}.
\end{proof}

With Lemmas 1 and 2, we are now able to derive the finite sample breakdown point of the PCS of $\pmb S^*$ and $\pmb t^*$.

\begin{prop47}
For $n>p+1>2$ and $\pmb X$ in general position, the finite sample breakdown point of $\pmb S^*$ is 
\begin{equation}
\varepsilon_n^*(\pmb S,\pmb X)=\frac{n-h+1}{n}. \nonumber
\label{eq:PCSbreakdown}
\end{equation}
\end{prop47}

\begin{proof} 
Consider first the situation where $c\leqslant n-h$. Then any $h$-subset $H^m$ of $\pmb X^{\varepsilon}$ contains 
at least $p+1$ members of $G$. In particular, for the chosen $h$-subset $H^*$, 
denote $G^*=\{H^*\cap \{1,\ldots,g\}\}$ with $|G^*|\geqslant p+1$. 
The members of $G^*$ are in general position  
so that $\ave_{i\in G^*}(\pmb x_i^{\varepsilon}-\pmb t)(\pmb x_i^{\varepsilon}-\pmb t)'\succ0$ 
 for any $\pmb t\in\mathbb{R}^p$. But 
$\pmb S^*(\pmb X^{\varepsilon})=\ave_{i\in H^*}(\pmb x_i^{\varepsilon}-\pmb t^*)(\pmb x_i^{\varepsilon}-\pmb t^*)'$ and 
$\ave_{i\in H^*\setminus G^*}(\pmb x_i^{\varepsilon}-\pmb t^*)(\pmb x_i^{\varepsilon}-\pmb t^*)'\succeq 0$ 
so that $\pmb S^*(\pmb X^{\varepsilon})\succ 0$ \citep[10.58]{ppcs:Seber} which implies that   $\underset{\pmb X^\varepsilon \in \mathcal{X^\varepsilon}}{\sup} \lambda_p(\pmb S^*(\pmb X^{\varepsilon}))>0$. 
Thus for breakdown to occur, the numerator of Equation~\eqref{eq:breakdownexp}, $\lambda_1(\pmb S^*(\pmb X^\varepsilon))$, must become unbounded.
Now, suppose that $\pmb S^*(\pmb X^\varepsilon)$ breaks down. This means that for any $U_4>\max_{i=1}^n||\pmb x_i||^2$, 
\begin{equation}
\underset{\pmb X^\varepsilon \in \mathcal{X^\varepsilon}}{\sup}\lambda_1(\pmb S^*(\pmb X^\varepsilon)) > U_4. 
\label{eq:contrabase}
\end{equation}
We will show that this leads to a contradiction. In Appendix 2 we show that
\begin{equation}
\lambda_1(\pmb S^*(\pmb X^\varepsilon)) \leqslant \underset{i \in H^*}{\max}||\pmb x^\varepsilon_i||^2.
\label{eq:append1}
\end{equation}
By Equations~\eqref{eq:contrabase} and~\eqref{eq:append1}, it follows that $\underset{\pmb X^{\varepsilon}\in\mathcal{X}^{\varepsilon}}{\sup} \underset{i \in H^*}{\max}||\pmb x^\varepsilon_i||^2 >U_4$.
Then, by Lemma 2 we have that $\underset{\pmb X^{\varepsilon}\in\mathcal{X}^{\varepsilon}}{\sup}I(H^*) > U_1/K$ with $K=\binom{h}{p}$, the number of all directions $\pmb a\in B(H^*)$.
In particular, this is also true for $U_1 > k(\pmb X)$, and by Lemma 1, $k(\pmb X) \geqslant I(H^g)$, implying  that $I(H^*)>I(H^g)$ $\forall H^g \in \mathcal{H}^g$, which is a contradiction to the definition of $H^*$. 
Since PCS is affine and shift equivariant, when $c > n-h$, we have by Equation~\eqref{eq:maxBP} that $\pmb S^{*}(\pmb X^{\varepsilon})$ breaks down.
\end{proof}

Equation~\eqref{eq:maxBP} and Theorem 1 show that the breakdown point of $\pmb S^*$ is maximal. The following theorem shows that the breakdown point of $\pmb t^*$ is also maximal.

\begin{prop48}
For $n>p+1>2$ and $\pmb X$ in general position, the finite sample breakdown point of  $\pmb t^*$ is 
\begin{equation}
\varepsilon_n^*(\pmb t,\pmb X)=\frac{n-h+1}{n}. \nonumber
\label{eq:meanbreakdown}
\end{equation}
\end{prop48}

\begin{proof} 
Consider first the situation where $c\leqslant n-h$. In Theorem 1, we showed that under this condition $\underset{i \in H^*}{\max}||\pmb x^\varepsilon_i|| < \infty$. 
Denote $\pmb x^\varepsilon_u$, where $u=\argmax_{i\in H^\ast} \|\pmb x_i^\varepsilon \|$.
Then, we have that $\pmb t^*(\pmb X^\varepsilon)$ does not break down since, by homogeneity of the norm and the triangle inequality,
\begin{equation}
||\pmb t^*(\pmb X^\varepsilon)|| \leqslant \frac{1}{h }\sum_{i\in H^\ast}\|\pmb x^\varepsilon_i \| \leqslant\frac{1}{h}\sum_{i=1}^h\|\pmb x^\varepsilon_u \| =\|\pmb x^\varepsilon_u\| <\infty. \nonumber
\end{equation}
For the case of $c > n-h$, Equation~\eqref{eq:maxBP} and affine equivariance imply that $\pmb t^*(\pmb X^\varepsilon)$ breaks down.
\end{proof}

We now relax the assumption that the members of $\pmb X^g$ lie in general position in $\mathbb{R}^{p}$ and substitute it by
the weaker condition that they all lie in general position on a common subspace in $\mathbb{R}^q$ for some $q<p$. Then PCS has the so-called exact fit property. 
Recall that $\pmb a$ are hyperplanes defined by $p$ points drawn from an $h$-subset $H\in\mathcal{H}$. If there are at least $h$ points lying on a subspace, then there exists an $h$-subset of points from this subspace.  
Let $\tilde{H}$ be this subset. Then, for any $\pmb a^+\in B(\tilde{H})$, both the numerator and denominator of Equation \eqref{pcs:crit1} equal zero and so $I(\tilde{H})=0$. Thus, we have without loss of generality that $H^*=\{i:d_{i}^2(\pmb a^+)=0\}$. In summary this means that if $h$ or more observations lie exactly on a subspace, the fit given by the observations in $H^*$ will coincide with this subspace, which is the defintion of the so-called 
 \textit{exact fit} property. Of course, since $|H^*|\geqslant h$, $H^*$ may contain outliers. Given $H^*$, one may proceed with the much simpler task of identifying the at most $|H^*|-h$ outliers in this smaller set of observations on a rank $q$ subspace spanned by the members of $H^*$.
\bigskip

\noindent\textbf{Appendix 1: Proof of affine equivariance of 
$(\pmb t^*(\pmb X^{\varepsilon}),\pmb S^*(\pmb X^{\varepsilon}))$}\\

\noindent Recall that a location vector $\pmb t(\pmb X^{\varepsilon})$ and a scatter matrix $\pmb S(\pmb X^{\varepsilon})$ are  affine equivariant if for any non-singular $p\times p$ matrices $\pmb B$ and $p$-vector $\pmb b$ it holds that:
\begin{eqnarray*}
\pmb t(\pmb B\pmb X^{\varepsilon}+\pmb 1_n\pmb b')&=&\pmb B\pmb t(\pmb X^{\varepsilon})+\pmb b\\ 
 \pmb S( \pmb B \pmb X^{\varepsilon}  +\pmb 1_n\pmb b')&=& \pmb B  \pmb S(\pmb X^{\varepsilon})  \pmb B'.
\end{eqnarray*}
Consider now affine transformations of $\pmb X^{\varepsilon}$:
\begin{eqnarray}\label{ae1}
\pmb y^{\varepsilon}_i=\pmb B\pmb x^{\varepsilon}_i+\pmb b',\quad i=1,\ldots,n.
\end{eqnarray}
for any non-singular $p\times p$ matrix $\pmb B$ and $p$-vector $\pmb b$. 
The directions $\pmb a_x$ ($\pmb a_y$) are orthogonal to hyperplanes 
 through $p$-subsets of $\pmb X^{\varepsilon}$ ($\pmb Y^{\varepsilon}$).
 Since $||\pmb x_i^{g}-\pmb x_j^{g}||>0\;\forall\;1\leqslant i<j\leqslant g$, 
 we can disregard all duplicated rows of $\pmb X^{\varepsilon}$ (and their partner duplicates in $\pmb Y^{\varepsilon}$), so that, w.l.o.g. all $p$-subsets of $\pmb X^{\varepsilon}$ 
 $(\pmb Y^{\varepsilon})$  yield a $p \times p$ matrix with unique rows. Let $p^0$ 
 be any such $p$-subset of $\{1:n\}$, and $\pmb a_x^0$ and $\pmb a_y^0$ the hyperplanes 
 through $\{\pmb x_{i}\}_{i\in p^0}$ and $\{\pmb y_{i}\}_{i\in p^0}$. Since Equation \eqref{ae1} 
 describes an affine transformation, it preserves collinearity:
 \begin{eqnarray}\label{ae2}
\{i:\pmb x_i'\pmb a^0_x=1\}=\{i:\pmb y_i'\pmb a^0_y=1\},
\end{eqnarray}
and the ratio of lengths of intervals on univariate projections~\citep[sec. 36]{hcs:W02}:
\begin{eqnarray}\label{ae3}
\frac{\displaystyle\ave_{i\in H^m}||\pmb x_i'\pmb a^0_x-\ave_{i\in p^0}\pmb x_i'\pmb a^0_x||}{\displaystyle\ave_{i\in H(\pmb a_x^0)}||\pmb x_i'\pmb a^0_x-\ave_{i\in p^0}\pmb x_i'\pmb a^0_x||}=\frac{\displaystyle\ave_{i\in H^m}||\pmb y_i'\pmb a^0_y-\ave_{i\in p^0}\pmb y_i'\pmb a^0_y||}{\displaystyle\ave_{i\in H(\pmb a_y^0)}||\pmb y_i'\pmb a^0_y-\ave_{i\in p^0}\pmb y_i'\pmb a^0_y||},
\end{eqnarray}
where for readability we denote  $H^{\pmb a}$ as $H(\pmb a)$. Equation \eqref{ae2} and~\eqref{ae3} imply
\begin{eqnarray}\label{ae4}
I(H^m,\pmb a_x^0)=I(H^m,\pmb a_y^0).
\end{eqnarray}
Equation \eqref{ae4} holds for any $p$-subset of $H^m$. Therefore, denoting 
$B_x(H^m)$ all directions perpendicular to hyperplanes through $p$ elements of
 $\{\pmb x^{\varepsilon}_i\}_{i\in H^m}$, and $B_y(H^m)$ the same but for $\{\pmb y^{\varepsilon}_i\}_{i\in H^m}$), it holds that
\begin{eqnarray*}
\ave_{\pmb a_x\in B_x(H^m)} I(H^m,\pmb a_x)=\ave_{\pmb a_y\in B_y(H^m)} I(H^m,\pmb a_y),\quad m=1\ldots,M
\end{eqnarray*}
and in particular for $H^m=H^*$. Since $\#\{H^*\}\geqslant p+1$, we have that if the members of 
$H^*$ lie in G.P. in $\mathbb{R}^p$, 
\begin{eqnarray}
\ave_{i\in H^*}( \pmb B \pmb x_i^{\varepsilon}+\pmb b')&=& \pmb B  \ave_{i\in H^*}( \pmb x_i^{\varepsilon})\pmb b, \nonumber\\
\cov_{i\in H^*}( \pmb B \pmb x_i^{\varepsilon}+\pmb b')&=& \pmb B  \cov_{i\in H^*}( \pmb x_i^{\varepsilon})  \pmb B'. \nonumber
\end{eqnarray}   
Hence, $(\pmb t^*(\pmb X^{\varepsilon}),\pmb S^*(\pmb X^{\varepsilon}))$ are affine equivariant. 
\medskip

\noindent\textbf{Appendix 2: Proof of Equation~\ref{eq:append1}}\\

\noindent Here, we show that $\lambda_1(\pmb S^*(\pmb X^\varepsilon)) \leqslant \underset{i \in H^*}{\max}||\pmb x^\varepsilon_i||^2$. 
The first eigenvalue of $\pmb S^*(\pmb X^{\varepsilon})$ is defined as 
$$\lambda_1(\pmb S^*(\pmb X^{\varepsilon}))=\cov_{i\in H^*}((\pmb x_i^{\varepsilon})'d)$$
for $d=\underset{||\tilde{d}||=1}{\argmax}\underset{i\in H^*}{\cov}((\pmb x_i^{\varepsilon})'\tilde{d})$.
Furthermore,
$$\cov_{i\in H^*}((\pmb x_i^{\varepsilon})'d)=\ave_{i\in H^*}(((\pmb x^{\varepsilon}_i)'d)^2)-(\ave_{i\in H^*}((\pmb x_i^{\varepsilon})'d))^2.$$
Hence, we have that
$$\cov_{i\in H^*}((\pmb x_i^{\varepsilon})'d)\leqslant\ave_{i\in H^*}(((\pmb x^{\varepsilon}_i)'d)^2).$$
We also have that
$$\ave_{i\in H^*}(((\pmb x^{\varepsilon}_i)'d)^2)\leqslant\max_{i\in H^*}(((\pmb x^{\varepsilon}_i)'d)^2)=\max_{i\in H^*}||(\pmb x^{\varepsilon}_i)'d||^2.$$
Using Cauchy-Schwartz,

$$\max_{i\in H^*}||(\pmb x^{\varepsilon}_i)'d||^2\leqslant(\max_{i\in H^*}||d||||\pmb x^{\varepsilon}_i)||)^2,$$
and $||d||=1$. Thus, $\lambda_1(\pmb S^*(\pmb X^{\varepsilon}))\leqslant\underset{i\in H^*}{\max}||\pmb x^{\varepsilon}_i||^2.$\\

\section{Acknowledgements}
\noindent The authors wish to acknowledge the helpful comments from three anonymous 
referees and the editor for improving this paper. 

\noindent Viktoria \"Ollerer would like to acknowledge the support of Research Fund KU Leuven GOA/12/014.
\bibliographystyle{model1-num-names}

\end{document}